\documentclass[12pt,oneside,final]{amsart}
\usepackage[cp1250]{inputenc}
\usepackage[T1]{fontenc}
\usepackage{amssymb}
\usepackage{color}
\renewcommand{\thefootnote}{}
\newtheorem{twr}{Twierdzenie}
\newtheorem{thr}[twr]{Theorem}
\newtheorem{prp}[twr]{Proposition}

\newtheorem{lm}[twr]{Lemma}

\renewcommand{\thefootnote}{}

\newcommand{\po}{{\partial\Omega}}

\newcommand{\cn}{\mathbb{C}^n}

\newcommand{\cj}{ \in\mathcal{C}^{1,1}}

\title[The  Monge-Amp\`{e}re equation]{The Monge-Amp\`{e}re equation on almost complex manifolds }
\author[S. Pliś]{Szymon Pliś}
\address{  Institute of Mathematics, Cracow University of Technology, Warszawska 24, 31-155
    Kraków, Poland
}

\email{splis@pk.edu.pl}

\begin{document}\thispagestyle{empty} \footnotetext{Partially supported by the project N N201 2683 35 of the Polish Ministry of
Science and Higher Education.
}\renewcommand{\thefootnote}{\arabic{footnote}}

\begin{abstract}
  We  study the Dirichlet problem for the
Monge-Amp\`ere equation on almost complex manifolds. We obtain the existence of the unique smooth solution of this problem in strictly pseudoconvex domains.
\end{abstract}

\maketitle

Let $(M,J)$ be an almost complex manifold of real dimension $2n$ (the definitions will be given in section \ref{sec: notion}).  N. Pali proved (in \cite{p}) that, as it is in the case of complex geometry, for plurisubharmonic functions the $(1,1)$ current $i\partial\bar\partial u$ is nonnegative\footnote{But it can be not closed!}. So for smooth plurisubharmonic function $u$ we have well defined Monge-Amp\`ere operator $(i\partial\bar\partial u)^n\geq0$ and we can study the complex Monge-Amp\`ere equation \begin{equation}\label{ma}
(i\partial\bar\partial u)^n=fdV,
\end{equation}  where $f\geq0$ and $dV$ is a (smooth) volume form 
.

Let $\Omega\Subset M$ be a strictly pseudoconvex domain of class $\mathcal{C}^{\infty}$. 
In this article we study the following Dirichlet problem for the Monge-Amp\`ere equation:
 \begin{equation}\label{DP}
\left\{
\begin{array}{l}
    u\in\mathcal{PSH}(\Omega)\cap \mathcal{C}^\infty(\bar\Omega)\\ 
    (i\partial\bar\partial u)^n=dV \;\mbox{ in }\;\Omega\\
    u=\varphi\;\mbox{ on }\;\partial\Omega
\end{array}
\right.
\end{equation}
where $\varphi\in\mathcal{C}^\infty(\bar\Omega)$. The main theorem is the following:

\begin{thr}\label{mine theorem}
There is a unique smooth plurisubharmonic solution $u$ of the problem (\ref{DP}).
\end{thr}

In \cite{c-k-n-s} the above theorem was proved for $\Omega\subset\cn$ with $J_{\rm st}$. Note that even in the integrable case it is not enough to assume that $\po$ is strictly pseudoconvex\footnote{We say that $\po$ is stictly pseudoconvex if is locally equal to $\{\rho=0\}$, for same smooth function $\rho$ such that $i\partial\bar{\partial}\rho>0$ and $\triangledown\rho\neq0$.}. Indeed, if $\Omega$ is the blow-up of a strictly pseudoconvex domain in $\cn$  in one point then $\po$ is strictly pseudoconvex.  But if $u\in\mathcal{PSH}(\Omega)\cap \mathcal{C}^\infty(\bar\Omega)$ the form $(i\partial\bar\partial u)^n$ is not a volume form.

 In case of $J$ not integrable McDuff construct a domain $\Omega$ with a non connected strictly pseudoconvex boundary (see \cite{m}).\footnote{If $\Omega$ is strictly pseudoconvex, then $\po$ is connected (see \cite{b-g}).} One can prove the above theorem for $\Omega$ not necessary strictly pseudoconvex, in almost the same way, if we assume that $\partial\Omega$ is strictly pseudoconvex, $\varphi$ is plurisubharmonic in neigborhood of $\bar{\Omega}$  and $dV\leq(i\partial\bar{\partial}\varphi)^n$ on $\Omega$. It is however not clear for the author, if there is an example of such $\varphi$ in McDuff's example (or in any other not strictly pseudoconvex domain with a strictly pseudoconvex boundary).

In the last section we will explain how Theorem \ref{mine theorem} gives the theorem of Harvey and Lawson about existing a continuous solution of The Dirichlet Problem for maximal functions. We even improve their result by proving that the solution is Lipschitz (if the boundary condition is enough regular).

\numberwithin{equation}{section}\numberwithin{twr}{section}
\section{notion}\label{sec: notion}
We say that $(M,J)$ is an almost complex manifold if $M$ is a manifold and $J$ is an ($\mathcal{C}^\infty$ smooth) endomorphism of the tangent bundle $TM$, such that $J^2=-\rm{id}$. The real dimension of $M$ is even in that case.

 We have then a direct sum decomposition $T_{\mathbb{C}}M=T^{1,0}M \oplus T^{0,1}M$, where $T_{\mathbb{C}}M$ is a complexification of $TM$, $$T^{1,0}M=\{X-iJX:X\in TM\}$$ and $$T^{0,1}M=\{X+iJX:X\in TM\}(=\{\zeta\in T_{\mathbb{C}}M:\bar\zeta\in T^{1,0}M\}).$$

Let $\mathcal{A}^k$ be the set of $k$-forms, i.e. the set of sections of $\bigwedge^k(T_{\mathbb{C}}M)^\star$ and let
$\mathcal{A}^{p,q}$ be the set of $(p,q)$-forms, i.e. the set of sections of $\bigwedge^p(T^{1,0}M)^\star\otimes_{(\mathbb{C})}~\bigwedge^q(T^{0,1}M)^\star$. Then we have a direct sum decomposition $\mathcal{A}^k=\bigoplus_{p+q=k}\mathcal{A}^{p,q}$. We denote the projections $\mathcal{A}^k\rightarrow\mathcal{A}^{p,q}$ by $\Pi^{p,q}$.

If $d:\mathcal{A}^k\rightarrow\mathcal{A}^{k+1}$ is (the $\mathbb{C}$-linear extension of) the exterior differential, then we define $\partial: \mathcal{A}^{p,q}\rightarrow\mathcal{A}^{p+1,q}$ as $\Pi^{p+1,q}\circ d$ and $\bar\partial: \mathcal{A}^{p,q}\rightarrow\mathcal{A}^{p,q+1}$ as $\Pi^{p,q+1}\circ d$.

We say that an almost complex structure $J$ is integrable, if satisfy any of the following (equivalent) conditions:\\i) $d=\partial+\bar\partial$;
\\ii) $\bar\partial^2=0$;
\\iii) $[\zeta,\xi]\in T^{0,1}M$ for vector fields $\zeta,\xi\in T^{0,1}M$.
\\ By the Newlander-Nirenberg Theorem $J$ is integrable if and only if it is induced by a complex structure.

In the paper $\zeta_1,\ldots,\zeta_n$ is always a (local) frame of $T^{1,0}$ . Let us put for a smooth function $u$
$$u_{p\bar q}=\zeta_p\bar{\zeta_q}u=u_{\bar qp}+[\zeta_p,\bar{\zeta_q}]u$$
 and
$$A_{p\bar q}=A_{p\bar q}(u)=u_{p\bar q}-[\zeta_p,\bar{\zeta_q}]^{0,1}u,$$
where $X^{0,1}=\Pi^{0,1}(X)$. Then for a smooth function $u$ we have (see \cite{p}):
 $$i\partial\bar\partial u=i\sum A_{p\bar q}d\zeta_p\wedge d\bar\zeta_q,$$
  where $d\zeta_1,\ldots, d\zeta_n,d\bar\zeta_1,\ldots, d\bar\zeta_n$ is a base of $(T_{\mathbb{C}}M)^\star$ dual to the base $\zeta_1,\ldots,\zeta_n,\bar{\zeta_1},\ldots,\bar{\zeta_n}$ of $T_{\mathbb{C}}M$.

Let $\mathbb{D}=\{z\in\mathbb{C}:|z|<1\}$. We say that a (smooth) function $\lambda:\mathbb{D}\rightarrow M$ is $J$-holomorphic or simpler holomorphic if $\lambda'(\frac{\partial}{\partial\bar z})\in T^{0,1}M$. The plenty of such disks show the following proposition from \cite{i-r}, where is stated for $C^{k',\alpha} $ class of $J$:

\begin{prp}\label{dyski}
 Let $0\in M\subset\mathbb{R}^{2n}$, $k,k'\geq1$. For $v_0,v_1,\ldots,v_k\in\mathbb{R}^{2n}$ close enough to 0, there is a holomorphic function $\lambda:\mathbb{D}\rightarrow M$, such that $\lambda(0)=v_0$ and $\frac{\partial^l\lambda}{\partial x^l}=v_l$, for $l=1,\ldots,k$. Moreover, we can choose $\lambda$ with $\mathcal{C}^1$ dependence on parameters $(v_0,\ldots,v_k)\in(\mathbb{R}^{2n})^{k+1}$, where for holomorphic functions  we consider $\mathcal{C}^{k'}$ norm.
\end{prp}

We can locally normalize coordinates with respect to a given holomorphic disc $\lambda$, that is we can assume that $\lambda(z)=(z,0)\in\cn$ and $J=J_{\rm st}$ on $\mathbb{C}\times\{0\}\subset\cn$, where $J_{\rm st}$ is the standard almost complex structure in $\cn$ (see section 1.2 in \cite{b}) and moreover we can assume that for every $J$-holomorphic $\mu$ such that $\mu(0)=0$ we have $\triangle \mu(0)=0$ (see \cite{r}).

An upper semi-continuous function $u$ on an open subset of $M$ is said to be $J$-plurisubharmonic or simpler plurisubharmonic if a function $u\circ\lambda$ is subharmonic for every holomorphic function $\lambda$. We denote the set of plurisubharmonic functions on $\Omega\subset M$ by $\mathcal{PSH}(\Omega)$. For a smooth function $u$ it means that a matrix $(A_{p\bar q})$ is nonnegative. Recently Harvey and Lawson proved that an upper semicontinuous function localy integrable $u$ is plurisubharmonic iff a current $i\partial\bar\partial u$ is nonnegative. We say that a function $u \in\mathcal{C}^{1,1}(\Omega)$ is strictly plurisubharmonic if for every $K\Subset\Omega$ there is $m>0$ such that $\omega\leq im\partial\bar\partial u$ a. e. in $K$, where $\omega$ is any hermitian metric\footnote{Hermitian metric is a smooth positive $(1,1)$ form.} on $\Omega$. If $u\in\mathcal{C}^2(\Omega)$ then the following conditions are equivalent:\\ i) $u$ is strictly plurisubharmonic;\\ 
ii) $i\partial\bar{\partial}u>0$;\\
iii) $u$ is plurisubharmonic and $(i\partial\bar{\partial}u)^n>0$.

We say that a domain $\Omega\Subset M$ is strictly pseudoconvex of class $\mathcal{C}^{\infty}$ (respectively of class $\mathcal{C}^{1,1}$) if there is a strictly plurisubharmonic function $\rho$ of class $\mathcal{C}^{\infty}$ (respectively of class $\mathcal{C}^{1,1}$) in a neighborhood of $\bar\Omega$,  such that $\Omega=\{\rho<0\}$ and $\triangledown\rho\neq0$ on $\partial\Omega$. In that case we say that $\rho$ is a defining function for $\Omega$. 

Let $z_0\in M$. The most basic example of a (strictly) plurisubharmonic function in neighborhood of $z_0$ is $u(z)={\rm dist}(z,z_0)$ (where ${\rm dist}$ is a distance in some rimmanian metric).  Domains $\Omega_\varepsilon=\{u<\varepsilon\}$ are strictly pseudoconvex of class  $\mathcal{C}^{\infty}$ for $\varepsilon>0$ small enough and they make a fundamental neighborhood system for $z_0$.

\section{ comparison principle}
In this section $\Omega\Subset M$ is a domain not necessary strictly pseudoconvex but such that there is a bounded function $\rho\in\mathcal{C}^2\cap\mathcal{ PSH}(\Omega)$. 

In pluripotential theory in $\cn$, the comparison principle is a very effective tool. We give here the most basic version for $J$-plurisubharmonic functions.

\begin{prp}[COMPARISON PRINCIPLE]\label{cp}
 If   $u$, $v\in\mathcal{C}^2(\bar\Omega)$ are such that $u$ is a plurisubharmonic functions, $(i\partial\bar{\partial} u)^n\geq(i\partial\bar{\partial} v)^n$ on the set $\{i\partial\bar{\partial} v>0\}$ and $u\leq v$ on $\po$, then $u\leq v$ in $\bar\Omega$.
\end{prp}

{\it Proof:} First, let us assume that $(i\partial\bar{\partial} u)^n>(i\partial\bar{\partial} v)^n$ on the set $\{i\partial\bar{\partial} v>0\}$ and a function $u-v$ takes its maximum in $z_0\in\Omega$. By proposition \ref{dyski} for small $\zeta\in T_{z_0}^{1,0}M $ there is a holomorphic disk $\lambda$ such that $\lambda(0)=z_0$ and $\lambda'(0)=\zeta  $. Hence at $z_0$ $$i\partial\bar\partial(v-u)(\zeta,\bar\zeta)=\frac{1}{4 }\bigtriangleup\left((v-u)\circ\lambda\right)(0)\geq 0$$ so we have $i\partial\bar\partial u\leq i\partial\bar\partial v$ and than we obtain $(i\partial\bar\partial u)^n\leq(i\partial\bar\partial v)^n$ which is the contradiction with our  first assumption.

In the general case we put $u'=u+\varepsilon(\rho-\sup_{\bar\Omega}\rho)$ and the lemma follows from the above case (with $u'$ instead of $u$).  $\;\Box$

In section \ref{smax} we will need a slight stronger version of the above proposition.

\begin{prp}\label{cp'}
Suppose that   $u$, $v\in\mathcal{C}^2(\bar\Omega)$ are such that $u$ is a plurisubharmonic functions, $(i\partial\bar{\partial} u)^n\geq(i\partial\bar{\partial} v)^n$ on the set $\{i\partial\bar{\partial} v>0\}$ and $u\leq v$ on $\po$. Then for any $H\in\mathcal{PSH}$, an inequality $$\varlimsup_{z\rightarrow z_0}(u+H-v)\leq0$$ for any $z_0\in\po$ implies $u+H\leq v$ on $\Omega$.
\end{prp}

{\it Proof:} Let $z_0\in\Omega$ be a point where a function $f=u+H-v$ attains a maximum and $\lambda$ is a holomorphic disk  such that $\lambda(0)=z_0$. Because  $H\circ\lambda$ is a subharmonic function one can find a sequence $t_k$ of nonzero complex numbers such that $$\lim_{k\rightarrow\infty}t_k=0$$ and $$4H\circ\lambda(0)\leq H\circ\lambda(t_k)+H\circ\lambda(it_k)+H\circ\lambda(-t_k)+H\circ\lambda(-it_k)$$ Hence
$$\bigtriangleup\left((v-u)\circ\lambda\right)(0)$$$$\geq\varlimsup_{k\rightarrow\infty}\frac{4f\circ\lambda(0)-f\circ\lambda(t_k)-f\circ\lambda(it_k)-f\circ\lambda(-t_k)-f\circ\lambda(-it_k)}{|t_k|^2}$$$$\geq 0.$$
From this we can obtain our result exactly as in the proof of the previous proposition. $\;\Box$

\section{a priori estimate}

In this section we will proof a $\mathcal{C}^{1,1}$ estimate for the smooth solution $u$ of the problem (\ref{DP}). By the general theory of elliptic equations (see for example \cite{c-k-n-s}) we obtain from this the $\mathcal{C}^{k,\alpha}$ estimate and then the existence of smooth solution. The uniqness follows from comparison principle.

Our proofs are close to \cite{c-k-n-s} but more complicated because of the noncommutativity of some vector fields.

\subsection{some technical preparation}

In this section we assume that $\Omega\Subset M$ is strictly pseudoconvex of class $\mathcal{C}^\infty$ with the defining function  $\rho$. Let us fix a hermitian metric $\omega$ on $M$. From now all norms, gradient and hessian are taken with respect to this metric or more precisely with respect to a rimannian metric which is given by $g(X,Y)=-\omega(X,JY)$ for vector fields $X$, $Y$.  

Let $f\in\mathcal{C}^\infty(\bar{\Omega})$ be such that $dV=f\omega^n$. Then locally our Monge-Amp\`ere equation $(i\partial\bar\partial u)^n=dV$ has a form:
$$ \det(A_{p\bar q})=gf,$$ where $g=\left(\det(\omega(\zeta_p,\bar\zeta_q))\right)^{-2}$. So if vectors $\zeta_1,\ldots,\zeta_n$ are ortonormal (i.e. $\omega(\zeta_p,\bar\zeta_q)=\delta_{pq}$) then $g=1$.

We will very often use the following elliptic operator
$$L=L_{\zeta}=A^{p\bar q}(\zeta_p\bar{\zeta_q}-[\zeta_p,\bar{\zeta_q}]^{0,1}).$$
Note that  for $X$, $Y$ vector fields we have
 \begin{equation*}\label{pp}
    X(\log gf)=A^{p\bar q}XA_{p\bar q},
 \end{equation*}
 \begin{equation*}\label{dp}
   XY(\log gf)=A^{p\bar q}XYA_{p\bar q}-A^{p\bar
    j}A^{i\bar q}(YA_{i\bar j})(XA_{p\bar q}),
 \end{equation*}
were $(A^{p\bar q})$ is the inverse of the matrix
$(\overline{A_{p\bar q}})$.



 In the Lemmas we will specify exactly how  {\it a priori estimates} depend on $\rho$, $f$ and $\varphi$. We should  emphasize that it depend strongly also on $M$, $J$ and $\omega$. The notion $C(A)$ really means that $C$ depend on an upper bound for $A$. $C$ always depends on    $m(\rho)$ what is define as the smallest constant $m>0$ such that $\omega\leq mi\partial\bar{\partial}\rho$ on $\Omega$.

In the proofs below $C$ is a constant under control, but it can change from a line to a next line.

\subsection{uniform estimate}

\begin{lm}\label{ogr}
We have $\|u\|_{L^\infty(\Omega)}\leq C$,\\ where
$\;C=C(\|\rho\|_{L^\infty(\Omega)},m(\rho),\|f\|_{L^\infty(\Omega)},\|\varphi\|_{L^\infty(\Omega)})$.
 \end{lm}
 \textit{Proof:} From the comparison principle and the maximum principle we have
 $$\|f^{1/n}\|_{L^\infty(\Omega)}m(\rho)\rho+\inf_\po\varphi\leq u\leq\sup_\po\varphi.\;\Box$$

\subsection{gradient estimate}

In two next lemmas we shall prove {\it a prori} estimate for first derivate.

\begin{lm}\label{gradientb}
We have $$\|u\|_{\mathcal{C}^{0,1}(\partial\Omega)}\leq C,$$ where
$\;C=C(\|\rho\|_{\mathcal{C}^{0,1}(\Omega)},m(\rho),\|f\|_{L^\infty(\Omega)},\|\varphi\|_{\mathcal{C}^{1,1}(\Omega)})$.
 \end{lm}


 \textit{Proof:} We can choose $A>0$ such that $Ai\partial\bar{\partial}\rho+i\partial\bar\partial\varphi\geq f^{1/n}\omega$ and $Ai\partial\bar{\partial}\rho\geq i\partial\bar\partial\varphi$. Thus by the comparison principle and the maximum principle we have $$\varphi+A\rho\leq u\leq \varphi-A\rho$$
 for $A$ large enough. So on the boundary we have $$|\nabla
 u|\leq|\nabla A\rho|+|\nabla \varphi|
.\;\Box$$

 \begin{lm}\label{gradient}
We have \begin{equation*}\label{grad}
\|u\|_{\mathcal{C}^{0,1}(\Omega)}\leq C,
\end{equation*}
 where
$\;C=C(\|\rho\|_{\mathcal{C}^{0,1}(\Omega)},m(\rho),\|f^{1/n}\|_{\mathcal{C}^{0,1}},\|u\|_{\mathcal{C}^{0,1}(\partial\Omega)})$.
 \end{lm}

\textit{Proof:}  Consider the function $v=\psi|\nabla u|^2$. We  assume that $v$ takes its maximum in $z_0\in\Omega$. We can choose $\zeta_1,\ldots,\zeta_n$, such that they are orthonormal in a neighbourhood of $z_0$, and the matrix $A_{p\bar q}$ is diagonal at $z_0$. From now on all formulas are assumed to hold at $z_0$.

 We have $Xv=0$ so that $X(|\nabla|^2)=-X\log\psi|\nabla u|^2$.
We can calculate
$$L(v)=L(\psi)|\nabla|^2+\psi L(|\nabla|^2)+A^{p\bar p}(\psi_p(|\nabla|^2)_{\bar p}+\psi_{\bar p}(|\nabla|^2)_p)$$
$$|\nabla|^2A^{p\bar p}(\psi_{p\bar p}-[\zeta_p,\bar{\zeta_p}]^{0,1}\psi-2\frac{|\psi_p|^2}{\psi})+\psi L(|\nabla|^2)$$
$$L(|\nabla|^2) =A^{p\bar p}((|\nabla|^2)_{p\bar p}-[\zeta_p,\bar{\zeta_p}]^{0,1}|\nabla|^2)$$
 $$= A^{p\bar p}\sum_k(u_{p\bar pk}u_{\bar k}+u_{k}u_{p\bar p\bar k}+|u_{pk}|^2+|u_{\bar pk}|^2 -[\zeta_p,\bar{\zeta_p}]^{0,1}u_ku_{\bar k}-u_k[\zeta_p,\bar{\zeta_p}]^{0,1}u_{\bar k})$$
$$A^{p\bar p}(u_{p\bar pk}-[\zeta_p,\bar{\zeta_p}]^{0,1}u_k)$$
 $$=A^{p\bar p}(u_{kp\bar p}-\zeta_k[\zeta_p,\bar{\zeta_p}]^{0,1}u+\zeta_p[\bar{\zeta_p},\zeta_k]u+[\zeta_p,\zeta_k]\bar{\zeta_p}u)$$ $$=(\log f)_k+A^{p\bar p}(\zeta_p[\bar{\zeta_p},\zeta_k]u+\bar{\zeta_p}[\zeta_p,\zeta_k]u+[[\zeta_p,\zeta_k],\bar{\zeta_p}]u-[[\zeta_p,\bar{\zeta_p}]^{0,1},\zeta_k]u)$$
then
we have
$$|A^{p\bar p}(u_{p\bar pk}-[\zeta_p,\bar{\zeta_p}]^{0,1}u_k)|$$ $$\leq C\left(\frac{\|f^{1/n}\|_{\mathcal{C}^{0,1}}}{f^{1/n}}+A^{p\bar p}\left(\sum_s(|u_{ps}|+|u_{p\bar s}|)+|\nabla u|\right)\right)$$
and similarly
$$|A^{p\bar p}(u_{p\bar p\bar k}-[\zeta_p,\bar{\zeta_p}]^{0,1}u_{\bar  k})|$$ $$\leq C\left(\frac{\|f^{1/n}\|_{\mathcal{C}^{0,1}}}{f^{1/n}}+A^{p\bar p}\left(\sum_s(|u_{ps}|+|u_{p\bar s}|)+|\nabla u|\right)\right)$$
so for the proper choice of $\psi$ (we can get $\psi=e^{A\phi}+B$ for $A$, $B$ large enough) we have $L(v)(0)>0$ and this is a contradiction with the maximality of $v$. $\;\;\;\Box$

\subsection{$\mathcal{C}^{1,1}$ estimate} 
Let us fix $P\in\po$. Now we will give $\mathcal{C}^{1,1}$ estimate in a point $P$ (which not depends on $P$). Estimate of $XYu(P)$ where $X,Y$ are tangent to $\po$ follows from the gradient estimate.

\begin{lm}\label{druganabrzegu}
Let $N\in T_PM$ be orthogonal to $\po$ such that $N\rho=-1$ and let $X$ be a vector field on a neighborhood of $P$ which is tangent to $\po$ on $\po$. We have \begin{equation*}\label{druganabrzeg}
|NXu(P)|\leq C,
\end{equation*}
 where 
$C=C(\|\rho\|_{\mathcal{C}^{0,1}(\Omega)},m(\rho),\|f^{1/n}\|_{\mathcal{C}^{0,1}},\|\varphi\|_{\mathcal{C}^{2,1}(\Omega)} , \|X\|_{\mathcal{C}^{0,1}},\|u\|_{\mathcal{C}^{0,1}(\Omega)})$.
 \end{lm}

\textit{Proof:} Let  $X_1,X_2,\ldots,X_{n}$ be (real) vector fields on $U$, a neighborhood of $P$, tangent at $P$ to $\po$, such that $X_1,JX_1,\ldots,X_n,JX_n$ is a frame. Consider the function $$v=X(u-\varphi)+B\rho+\sum_{k=1}^n|X_k (u-\varphi)|^2-A({\rm dist}(P,\cdot))^2.$$
Let $V\Subset U$ be a neighborhood of $P$ and  $S=V\cap\Omega$. For $A$ large enough $v\leq0$ on $\partial S$.

Our goal is to show that for $B$ large enough we have $v\leq0$ on $\bar S$. Let $z_0\in S$ be a point where $v$ attains a maximum and let $\zeta_1,\ldots,\zeta_n$ be orthonormal and such that $(A_{p\bar q})$ is diagonal. From now on all formulas are assumed to hold at $z_0$. Let us calculate:
$$ m(\rho)L(\rho)\geq \sum A^{p\bar p}$$and
$$L(-X\varphi-A({\rm dist}(P,\cdot))^2)\geq-C\sum A^{p\bar p},$$ hence for $B$ large enough
$$L(B\rho-X\varphi-A({\rm dist}(P,\cdot))^2)\geq \frac{B}{2m(\rho)}\sum A^{p\bar p}.$$
To estimate $L(Xu+\sum_{k=1}^{n}|X_k (u-\varphi)|^2)$ let us first consider $Y\in\{X,X_1,\ldots X_{n}\}$ and  calculate 
$$L(Yu)=A^{p\bar q}(\zeta_p\bar\zeta_q Yu-[\zeta_p,\bar\zeta_q]^{0,1}Yu)$$ $$=Y\log f +A^{p\bar q}(\zeta_p[\bar\zeta_q, Y]u+[\zeta_p, Y]\bar\zeta_qu-[ [\zeta_p,\bar\zeta_q]^{0,1},Y]u).$$
There are $\alpha_{q,k},\beta_{q,k}\in\mathbb{C}$ such that
$$[\bar\zeta_q, Y]=\sum_{k=1}^{n}\alpha_{q,k}\bar\zeta_k+\beta_{q,k}X_k$$
and so
$$A^{p\bar q}\zeta_p[\bar\zeta_q, Y]u=\sum_{q}\alpha_{q,q}+\sum_{k=1}^{n}A^{p\bar p}\beta_{p,k}\zeta_pX_ku+A^{p\bar p}Z_p u$$ where $Z_p$ are  vectors field  which are under control. This gives us
$$|A^{p\bar q}\zeta_p[\bar\zeta_q, Y]u|\leq C A^{p\bar p}(1+\sum_{k}|\beta_{p,k}\zeta_pX_ku|).$$
 In a similar way we can estimate $A^{p\bar q}[\zeta_p, Y]\bar\zeta_qu$ and we obtain
$$|L(Yu)|\leq CA^{p\bar p} (1+\sum_k|\zeta_pX_ku|).$$
Therefore
$$L(Xu+\sum_k|X_k (u-\varphi)|^2)$$
   $$\geq A^{p\bar p}\sum_{k=1}^{n}(\zeta_pX_k (u-\varphi))(\bar\zeta_{ p}X_k (u-\varphi))-CA^{p\bar p} (1+\sum_k|\zeta_pX_ku|)$$
  $$\geq A^{p\bar p}\sum_k|\zeta_pX_k u|^2-CA^{p\bar p} (1+\sum_k|\zeta_pX_ku|).$$

Now for $B$ large enough, since $L(v)(z_0)>0$,  we have  contradiction with maximality of $v$. Hence $v\leq0$ on $S$ and so $NXu(P)\leq C$ $\;\Box$.

\begin{lm}\label{druganormalnanabrzegu}
Let $X$ be a vector field orthogonal to $\po$ at $P$. We have \begin{equation*}\label{druganabrze}
\|XXu(P)\|\leq C,
\end{equation*}
 where  
$$C=C(\|\rho\|_{\mathcal{C}^{2,1}(\Omega)},m(\rho),\|f^{1/n}\|_{\mathcal{C}^{0,1}},\|f^{-1}\|_{L^\infty({\Omega})},\|\varphi\|_{\mathcal{C}^{3,1}(\Omega)} , \|X\|_{\mathcal{C}^{0,1}},\|u\|_{\mathcal{C}^{0,1}(\Omega)}).$$
 \end{lm}

\textit{Proof:} By the previous Lemma it is enough to prove that \begin{equation*}
|\zeta|^2\leq C\left(\zeta\bar\zeta -[\zeta,\bar\zeta]^{0,1}\right)u(P)
\end{equation*} for every vector field $\zeta\in T^{1,0}M$ tangent (at $P$) to $\po$.

Because
our argue will be local we can assume that $P=0\in\cn$. Let $\zeta_1,\zeta_2,\ldots\zeta_n\in T^{1,0}$ be a orthonormal frame in a neigborhood of 0 such that $\zeta_k\rho=-\delta_{kn}$. We can assume that $\zeta_1=\zeta$. By the strictly pseudoconvexity we have $(\zeta\bar\zeta -[\zeta,\bar\zeta]^{0,1})\rho(P)\neq0,$ so we can also assume that $(\zeta\bar\zeta -[\zeta,\bar\zeta]^{0,1})\varphi(P)=0$. 

From the strictly pseudoconvexity and using the proposition \ref{dyski} (for $k=2$) we can choose $J$-holomorphic disk $\lambda$ such that
$\lambda(0)=0$, $\frac{\partial \lambda}{\partial z}(0)=a\zeta$  and \begin{equation}\label{sty}
\rho\circ\lambda(z)= b|z|^2+O(|z|^3)
\end{equation} for some $a,b>0$. In particular we have
\begin{equation}\label{styk}
|z|^2\leq C{\rm dist}(\lambda(z),\bar\Omega).
\end{equation} Indeed, for $a>0$ smal enoughl, by the proposition \ref{dyski} (for $k=1$) there is a $J$-holomorphic disk $\tilde{\lambda}$ such that $\tilde\lambda(0)=0$ and $\frac{\partial \tilde{\lambda}}{\partial z}(0)=a\zeta$. Then (by changing coordinates) we can assume $\zeta J(0)=0$, so for any $J$-holomorphic disk $\lambda$ such that $\lambda(0)= 0$ and $\frac{\partial {\lambda}}{\partial z}(0)=a\zeta$ we have $\frac{\partial^2\lambda}{\partial x^2}(0)=-J\frac{\partial^2\lambda}{\partial x\partial y}(0)=-\frac{\partial^2\lambda}{\partial y^2}(0)$. Now if we put $$\frac{\partial^2\lambda}{\partial x^2}(0)=a^2(0,4\frac{\partial^2\rho}{\partial x_1^2}(0)-2(\zeta\bar\zeta -[\zeta,\bar\zeta]^{0,1})\rho(0),-4\frac{\partial^2\rho}{\partial x_1\partial x_2}(0))$$ we obtain (\ref{sty}) with $b=a^2(\zeta\bar\zeta -[\zeta,\bar\zeta]^{0,1})\rho(0)$.

 Once again changing coordinates we may assume  $\lambda(z_1)=(z_1,0)$, $\zeta_k(0)=\frac{\partial}{\partial z_k}$ for $k=1,\ldots,n$ and for every $J$-holomorphic disk $\mu$ such that $\mu(0)=0$ we have  
\begin{equation}\label{DSTmormalisation}
 \frac{\partial^2\mu}{\partial  z\partial\bar z}(0)=0.
\end{equation}

 We can find  a holomorphic cubic polynomial $p_1$ and a complex number $\alpha$ such that 
$$\varphi(z)=\varphi(0)+\varphi'(0)(z)$$  $$+{\rm Re}\left(\sum_{p=2  }^n\frac{\partial^2\varphi}{\partial z_1\partial\bar z_p}z_1\bar z_p+p_1(z)+\alpha z_1|z_1|^2\right)+O(|z_1|^4+|z_2|^2+\ldots+|z_n|^2).$$
By (\ref{styk}) we can choose another cubic polynomial $ p_2$ and numbers $\beta_1,\beta_2,\ldots,\beta_n\in\mathbb{C}$, $\beta_1>0$ such that  $${\rm Re}z_n={\rm Re}\left(\sum_{p=1}^n\beta_pz_1\bar z_p+ p_2(z)\right) +O(|z_1|^3+|z_2|^2+\ldots+|z_n|^2) \hbox{ on } \partial\Omega.$$   Then we obtain

$$\varphi(z)+\varphi(0)$$
$$=\varphi'(0)(z)+{\rm Re}\left(\sum_{p=2}^na_pz_1\bar z_p+ p_3(z)\right)+O(|z_2|^2+\ldots+|z_n|^2) $$
for some numbers $a_2,\ldots,a_n\in\mathbb{C}$ and a new cubic polynomial $p_3$, hence \begin{equation}
\label{unabrzegu}
u(z)-u(0)={\rm Re}\left(\sum_{p=2}^na_pz_1\bar z_p+ p_4(z)\right)+O(|z_2|^2+\ldots+|z_n|^2)\end{equation} for $z\in\partial\Omega$ and same polynomial $p_4$.

Let $B>1$ \footnote{Constants $C$ below not depend on this constant $B$.} and $D=B^{-1}\max\{|a_2|,\ldots,|a_n|\}$. By the proposition \ref{dyski} 
 (again for $k=2$) there is a family of $J$-holomorphic disks $g_w:\mathbb{D}\rightarrow\cn$, $w\in\mathbb{C}^{n-1}$ such that
 $$g_w(0)=(0,w),$$
  $$\frac{\partial g_w}{\partial z}(0)=(1,-\frac{ a_2}{B},\ldots,-\frac{ a_n}{B}),$$ 
 $$\|g_w-\lambda\|_{\mathcal{C}^4}\leq C(|w|+D)$$
and a function $G:\mathbb{D}\times\mathbb{C}^{n-1}\rightarrow\cn$ given by $G(z,w)=g_w(z)$ is of class $\mathcal{C}^4$. Than we have \begin{equation}\label{dysk}
|g_w(z)-(w_1+z,w_2-\frac{a_2z}{B},\ldots,w_n-\frac{a_nz}{B})|<C|z|^2(|w|+D)\end{equation} for  $z\in\mathbb{D}$ and  \begin{equation}\label{dysk'}|z|<C(\sqrt{|w|}+D)\end{equation} if  $g_w(z)\in\Omega.$

We can choose domains $ U\subset\mathbb{D}$, $V\subset\mathbb{C}^{n-1}$, $W\subset\mathbb{C}^{n}$ such that $0\in W$, $G(\partial U\times V)\cap\bar\Omega=\emptyset$ and $G|_{U\times V}$ is a diffeomorphism onto $W$.

Let $$h(g_w(z))={\rm Re}p_w(z)+AB |w|^2+\varepsilon\rho$$ where $A,\varepsilon>0$ and $p_w$ is a holomorphic cubic polynomial in one variable such that $${\rm Re} p_4(g_w(z))={\rm Re}p_w(z)+a_w|z|^2+{\rm Re}b_wz|z|^2+O(|z|^4)$$ for some $a_w\in\mathbb{R}$, $b_w\in\mathbb{C}$.
Note that $|b_w|< C(|w|+D)$ and by (\ref{DSTmormalisation}) $|a_w|<C|w|$. Thus enlarging $A$ (if necessary) and using also (\ref{dysk'}) we obtain 
\begin{equation}\label{p4}
{\rm Re}p_4(g_w(z))\leq h(g_w(z))+\frac{1}{2}D^2|z|^2 
\end{equation} 
on $\partial\Omega$.

 By inequalities (\ref{dysk}) and (\ref{dysk'})  we have
$$2\sum_{k=2}^n{\rm Re}a_kg_w^1(z)g_w^k(z)=\sum_{k=2}^nB(|-\frac{a_kg_w^1(z)}{B}-g_w^k(z)|^2-|\frac{a_kg_w^1(z)}{B}|^2-|g_w^k(z)|^2)$$
$$\leq B(|w|^2-D^2|z|^2+C|z|^4(|w|^2+D^2))\leq CB|w|^2-\frac{1}{2}D^2|z|^2$$ for $B$ enough large.
 By an above estimate (\ref{unabrzegu}),  and (\ref{p4}) we obtain that, if $A$ is enough large then $h\geq u-u(0)$ on $\partial\Omega\cap W$.  Enlarging $A$ again we can assume $h\geq u-u(0)$ on $\partial S$  where $S=\Omega\cap W$. Since $i\partial\bar\partial h$ is under control for $\varepsilon$ enough small we get an inequality $$(i\partial\bar\partial h)^n<(i\partial\bar\partial u)^n$$ on the set $S\cap\{i\partial\bar\partial h>0\}$. This by the Comparison Principle gives us $ h\geq u-u(0)$ on $S$.  
Note that $h_N\geq u_N$, $\varphi_{1\bar1}=0$ and $\varphi_N=h_N-\varepsilon\rho_N=h_N+\varepsilon$, so  we can conclude that $$u_{1\bar1}=u_{1\bar 1}-\varphi_{1\bar1}=(\varphi_N-u_N)\rho_{1\bar1}\geq\varepsilon\rho_{1\bar1}.  \;\;\;\Box$$

Finally we will obtain the interior $\mathcal{C}^{1,1}$ estimate, which together with previous lemmas gives us a full $\mathcal{C}^{1,1}$ estimate. By a standard argumentation this end the proof of Theorem 1.
\begin{lm}\label{drugawewnetrzu}
We have \begin{equation*}\label{drugawewnetrzu}
\|Hu\|_{L^\infty(\Omega)}\leq C,
\end{equation*}
 where
$Hu$ is a Hessian of $u$ and $$C=C(\|\rho\|_{\mathcal{C}^{0,1}(\Omega)},m(\rho),\|f^{\frac{1}{2n}}\|_{\mathcal{C}^{1,1}},\|u\|_{\mathcal{C}^{0,1}(\Omega)},\|Hu\|_{L^\infty(\po)}).$$
 \end{lm}

\textit{Proof:}
Let us define $M$ as the biggest eigenvalue of the Hessian $Hu$. We will show that the function
$$\Lambda=\psi e^{K|\nabla u|^2}M,$$
where $K^{-1}$ is large enough, does not attain maximum in $\Omega$.

Assume that a maximum of the function $\Lambda$ is attained at a point $z_0\in~\Omega$ (otherwise we are done).  There are $\zeta_1,\ldots,\zeta_n\in T_{z_0}^{1,0}M$ orthonormal at $z_0$ such that the matrix $(A_{p\bar q})$ is diagonal at $z_0$. Let $X\in T_{z_0}M$ be such that $\|X\|=1$ and $M=H(X,X)$. We can normalize coordinates near $z_0$ such that $z_0=0\in\cn$, $X=\frac{\partial}{\partial x_1}(0)$ and $J(z,0)=J_{st}$ for small $z\in\mathbb{C}$. Let us extend $X$ as $\frac{\partial}{\partial x_1}$ and  then in a natural way we  can extend $\zeta_1,\ldots,\zeta_n$ to some neigberhood $U$ of $0$ such that $[\zeta_k,X](0)=0$ and $[\zeta_k,\zeta_k](0)=0$ for $k=1,\ldots,n$. Indeed, first on $U\cap\mathbb{C}\times\{0\}$ we can put $\zeta_k$  as the same linear combination of vectors $\frac{\partial}{\partial z_1}\ldots,\frac{\partial}{\partial z_n}$ as in $0$. Then for some small $a>0$ we can take (for  $\zeta_k$ not tangent to $\mathbb{C}\times\{0\}$) $J$-holomorphic disks $d_k:\mathbb{D}\rightarrow U$ such that $d_k(0)=0$ and $\frac{\partial d_k}{\partial z}(0)=a\zeta_k(0)$, and on the image of $d_k$ we can put $\zeta_k(w)=a^{-1}\frac{\partial d_k}{\partial z}(d_k^{-1}(w))$. On the end we extend the vector fields on whole $U$.

 Let $$v=\psi e^{K|\nabla u|^2}\frac{Hu(\frac{\partial}{\partial x_1},\frac{\partial}{\partial x_1})}{|\frac{\partial}{\partial x_1}|^2}=\Psi e^{K|\nabla u|^2}(u_{x_1x_1}+Tu) \hbox{ on } U,$$ where $\Psi=\frac{\psi}{|\frac{\partial}{\partial x_1}|^2}$ and $T$ is a vector field (which is under control), then also a function $v$ has a maximum at $0$. Let us put $\mu= u_{x_1x_1}+Tu$ (then $|\frac{\partial}{\partial x_1}(0)|^2\mu(0)=M(0)$). Note that $XYu\leq C\mu$ for vector fields $X$, $Y$ (which are under control). 
  Assume $\mu>1$ (otherwise we have $\Lambda<C$, so we are done).

From now all formulas are assumed to hold at $0$. We will estimate $L(v)$ from below:
$$ L(v)= L(\Psi e^{K|\nabla u|^2})\mu+ \Psi e^{K|\nabla u|^2}L(\mu) -2A^{p\bar p}\frac{(\Psi e^{K|\nabla u|^2})_p(\Psi e^{K|\nabla u|^2})_{\bar p}}{\Psi e^{K|\nabla u|^2}}$$
To estimate the first term let us calculate
$$L(\Psi e^{K|\nabla u|^2})$$ $$=e^{K|\nabla u|^2}A^{p\bar p}(\Psi_{p\bar p}+2K{\rm Re}(\Psi_p(|\nabla u|^2)_{\bar p})+K\Psi(|\nabla u|^2)_{p\bar p}+K^2\Psi|(|\nabla u|^2)_p|^2),$$
$$A^{p\bar p}(|\nabla u|^2)_{p\bar p}$$ $$=A^{p\bar p}\sum_k((\zeta_p\bar\zeta_{ p }\eta_ku)u_{\bar k}+u_{ k}(\zeta_p\bar\zeta_{ p }\bar\eta_ku)+(\bar\zeta_p\eta_ku)(\zeta_p\bar\eta_ku)+(\zeta_p\eta_ku)(\bar\zeta_p\bar\eta_ku))$$
$$=\sum_k((\log\tilde f)_ku_{\bar k}+(\log\tilde  f)_{\bar k}u_k)$$ $$+A^{p\bar p}\sum_k((\zeta_p[\zeta_{\bar p },\eta_k]u)u_{\bar k}+([\zeta_p,\eta_k]u_{\bar p })u_{\bar k}+(\zeta_p[\bar\zeta_{ p },\bar\eta_k]u)u_{ k}+([\zeta_p,\bar \eta_k]u_{\bar p })u_{k}+(\eta_k+\bar\eta_k)[\zeta_p,\bar\zeta_p]^{0,1}u)$$ $$+A^{p\bar p}\sum_k((\bar\zeta_p\eta_k u)( \zeta_p\bar\eta_k u)+(\zeta_p\eta_k u)(\bar\zeta_p\bar\eta_k u)),$$
where $\eta_1,\ldots\eta_n$ is an orthonormal frame such that $\eta_k(0)=\zeta_k(0)$. Therefore we have
$$A^{p\bar p}(|\nabla u|^2)_{p\bar p}\geq -C+A^{p\bar p}\frac{1}{2}\sum_k((\bar\zeta_p\zeta_k u)( \zeta_p\bar\zeta_k u)+(\zeta_p\zeta_k u)(\bar\zeta_p\bar\zeta_k u)-C),$$
hence
$$L(\Psi e^{K|\nabla u|^2})\geq A^{p\bar p}\left(\Psi_{p\bar p}-CK(|\frac{\Psi_p}{\Psi}|^2+1)+ \frac{1}{2}\Psi K\sum_k((\bar\zeta_p\zeta_k u)( \zeta_p\bar\zeta_k u)+(\zeta_p\zeta_k u)(\bar\zeta_p\bar \zeta_k u))\right). $$
Let us start the calculation for the second term
$$L(\mu)=L(u_{x_1x_1})+L(Tu),$$
$$L(Tu)\leq T(\log f)-C(\mu+1)\sum A^{p\bar p},$$
$$L(u_{x_1x_1})= (\log f)_{x_1x_1}+A^{p\bar p}A^{q\bar q}|X(\zeta_p\bar \zeta_q-[\zeta_p,\bar \zeta_q]^{0,1}) u|^2$$ $$+A^{p\bar p}(\zeta_p[\bar\zeta_p,X]Xu+[\zeta_p,X]\bar\zeta_pXu+X\zeta_p[\bar\zeta_p,X]u+X[\zeta_p,X]\bar\zeta_pu+XX[\zeta_p,\bar\zeta_p]^{0,1}u)$$
$$= (\log f)_{x_1x_1}+A^{p\bar p}A^{q\bar q}|X(\zeta_p\bar \zeta_q-[\zeta_p,\bar \zeta_q]^{0,1}) u|^2$$ $$+A^{p\bar p}([\zeta_p,[\bar\zeta_p,X]]Xu+X[\zeta_p,[\bar\zeta_p,X]]u+[X,[\bar\zeta_p,X]]\zeta_pu+[X,[\zeta_p,X]]\bar\zeta_pu) $$ $$+A^{p\bar p}(X[X,[\zeta_p,\bar\zeta_p]^{0,1}]u+[X,[\zeta_p,\bar\zeta_p]^{0,1}]Xu)$$
$$\geq-C(\mu+1)\sum A^{p\bar p}$$
and we obtain
$$L(\mu)\geq-C\mu\sum A^{p\bar p}.$$
Now we come to the last term
$$-2A^{p\bar p}\frac{(\Psi e^{K|\nabla u|^2})_p(\Psi e^{K|\nabla u|^2})_{\bar p}}{\Psi e^{K|\nabla u|^2}}$$
$$=-2A^{p\bar p}e^{K|\nabla u|^2}(\frac{\Psi_p\Psi_{\bar p}}{\Psi}+K\Psi_p(|\nabla u|^2)_{\bar p}+K\Psi_{\bar p}(|\nabla u|^2)_p+K^2\Psi(|\nabla u|^2)_{\bar p}(|\nabla u|^2)_{ p})$$
$$\geq -2e^{K|\nabla u|^2}A^{p\bar p}\left(\frac{|\Psi_p|^2}{\Psi}+CK\eta|\Psi_p|+K^2\sum_k((\bar\zeta_p\zeta_k u)( \zeta_p\bar\zeta_k u)+(\zeta_p\zeta_k u)(\bar\zeta_p\bar \zeta_k u)))\right).$$

Therefore we can write $$ L(v)\geq \eta e^{K|\nabla u|^2}A^{p\bar p}(\Psi_{p\bar p}-\frac{|\Psi_p|^2}{\Psi}-CK(|\Psi_p|+|\frac{\Psi_p}{\Psi}|^2+1))$$
$$+ e^{K|\nabla u|^2}A^{p\bar p}(\frac{K\mu}{2}-2K^2)\sum_k((\bar\zeta_p\zeta_k u)( \zeta_p\bar\zeta_k u)+(\zeta_p\zeta_k u)(\bar\zeta_p\bar \zeta_k u))).$$ For a proper choice of $K$ and $\psi$ ($K<\frac{1}{4}$ and $\psi=e^{A\phi}+B$ where $\phi>1$ is strictly plurisubharmonic in a neigborhood of $\bar{\Omega}$ and $A$, $B$ are large enough) we can conclude $L(v)>0$ and it is a contradiction with the maximality of $v$. $\;\Box$



\section{ maximal plurisubharmonic functions}\label{smax}

We say that a function $u\in\mathcal{PSH}(\Omega)$ is maximal if for every function $v\in\mathcal{PSH}(\Omega)$ such that $v\leq u$
outside a compact subset of $\Omega$ we have $v\leq u$ in $\Omega$. 

Now we want to find the solution to the following Dirichlet problem: 
\begin{equation}\label{MP}
\left\{
\begin{array}{l}
    u\in\mathcal{PSH}(\Omega)\cap \mathcal{C}  (\bar\Omega)\\ 
    u \mbox{ is maximal }\\
    u=\varphi\;\mbox{ on }\;\partial\Omega
\end{array}
\right.
\end{equation}
where $\Omega$ is a strictly pseudoconvex domain of class $\mathcal{C}^{1,1}$ and $\varphi$ is a continuous function on $\partial\Omega$.

\begin{prp}\label{Lip}
If $\varphi\in\mathcal{C}^{1,1}(\bar{\Omega})$ then there is a unique solution $u\in\mathcal{C}^{0,1}(\bar{\Omega})$ of the problem (\ref{MP}) and $$\|u\|_{\mathcal{C}^{0,1}(\bar{\Omega})}\leq C=C(\|\rho\|_{\mathcal{C}^{0,1}(\Omega)},m(\rho),\|\varphi\|_{\cj}).$$
\end{prp}

{\it Proof:} The uniqueness is a consequence of the definition.

 To prove the existence first we assume that $\rho$ is smooth. There are an increasing sequence $\varphi_k$ of smooth functions such that $\varphi_k$ tends to $\varphi$ in $\mathcal{C}^{1,1}$ norm. By Theorem \ref{mine theorem} there is a solution $u_k$ of the following Dirichlet Problem
$$\left\{
\begin{array}{l}
    u_k\in\mathcal{PSH}(\Omega)\cap \mathcal{C}^\infty(\bar\Omega)\\ 
    (i\partial\bar\partial u_k)^n=\frac{1}{k^n}(i\partial\bar\partial \rho)^n \;\mbox{ in }\;\Omega\\
    u_k=\varphi_k\;\mbox{ on }\;\partial\Omega\;.
\end{array}
\right.$$ By Lemma \ref{gradient}  $\|u_k\|_{\mathcal{C}^{0,1}(\bar{\Omega})}\leq C(\|\rho\|_{\mathcal{C}^{0,1}(\Omega)},m(\rho),\|\varphi\|_{\cj})$. Now we can put $$u:=\lim_{k\rightarrow\infty} u_k.$$ 
It is enough to show that $u$ is a maximal function. Let a function $v\in\mathcal{PSH}(\Omega)$ be smaller than $u$ outside a compact subset of $\Omega$. From the comparison principle (Proposition \ref{cp'}) we obtain $$v+\frac{\rho}{k}-\sup_{\po}(\varphi-\varphi_k)\leq u_p$$
for $p\geq k$. Taking the limit we conclude that $v\leq u$ in $\Omega$.

In the general case we can assume that $\varphi$ is a plurisubharmonic function on $\Omega$ (by adding $A\rho$ for $A$ enough large). We can aproximate $\Omega$ by an increasing sequence of smooth strictly pseudoconvex domains $\Omega_k$ such that $\bigcup_k\Omega_k=\Omega $ and $\|\rho\|_{\mathcal{C}^{0,1}(\Omega)},m(\rho)$ are under control where $\rho_k$ are strictly plurisubharmonic smooth defining functions for $\Omega_k$. Let $u_k$ be a solution of the following Dirichlet Problem
$$\left\{
\begin{array}{l}
    u_k\in\mathcal{PSH}(\Omega_k)\cap \mathcal{C}  (\bar\Omega_k)\\ 
    u_k \mbox{ is maximal }\\
    u_k=\varphi\;\mbox{ on }\;\partial\Omega_k\;.
\end{array}
\right.$$ Then $u_k\geq\varphi$ hence it is an increasing sequence and again we can put $$u:=\lim_{k\rightarrow\infty} u_k.$$ If $v$ is as above, for every $\varepsilon>0$ we have $v-\varepsilon\leq u_k$ outside a compact set for $k$ large enough. So we obtain $v\leq u$ and we conclude that $u$ is a maximal function as in the statement. $\;\Box$

Note that in the above proposition it is not enough to assume that $\varphi$ is $\mathcal{C}^{1,\alpha}$ regular. Indeed, one can show that if $\Omega$ is strictly pseudoconvex, $P\in\po$ and $\varphi(z)\leq\varphi(P)-({\rm dist}(z,p))^{1+\alpha}$, then a solution of (\ref{MP}) is not H\"{o}lder continuous with the exponent greater than $\frac{1+\alpha}{2}$.

\begin{thr}[Harvey, Lawson \cite{h-l}]\label{maksymalne}
There is a unique solution $u$ of the problem (\ref{MP}).
\end{thr}

{\it Proof:} Let $\varphi_k $ be an increasing sequence of smooth functions on $\bar{\Omega}$ such that $\lim_{k\rightarrow\infty}\varphi_k=\varphi$. By Proposition \ref{Lip} there is a sequence $u_k$ of  solutions of (\ref{MP}) with boundary conditions $\varphi_k$ (instead of $\varphi$). Because $$u_k\leq u_p\leq u_k+\sup_{\po}(\varphi-\varphi_k)$$ for $p\geq k$, the sequence $u_k$ is a Cauchy sequence in $\mathcal{C}(\bar{\Omega})$. Similar as in the previous proof we can conclude that its limit $u$ is a solution of the problem (\ref{MP}). $\;\Box$ 

Note that we can also prove the above theorem directly from Theorem \ref{mine theorem}.

The following proposition shows that being a continuous maximal plurisubharmonic function is a local property.

\begin{prp}\label{lok}
Let $\Omega\subset M$ and $u\in\mathcal{PSH}(\Omega)$.    Then\\ i) If $u$ is maximal then $u|_U$ is maximal for every $U\subset\Omega$;\\ ii) If $\Omega$ is such that there is a bounded strictly plurisubharmonic function $\rho\in\mathcal{C}^2(\Omega)$, $u$ is continuous and  every point in $\Omega$ has a neighborhood $U$ such that $u|_U$ is maximal, then $u$ is maximal. 
\end{prp}

{\it Proof:} i) Suppose that $v\in\mathcal{PSH}(U)$ is such that $v\leq u$ outside a compact subset of $U$. Then $\max\{u,v\}\in\mathcal{PSH}(\Omega)$ and we obtain $v\leq\max\{u,v\}\leq u$ on $U$.

ii)  We can assume that $\rho<0$. Let $\varepsilon>0$, $v\in\mathcal{PSH}(\Omega)$ and let $z_0\in\Omega$ be a point where a function $v+\varepsilon\rho-u$ attains its maximum. By i) there is a stictly pseudoconvex domain $\tilde{\Omega}\subset\Omega$ with a smooth plurisubharmonic defining function $\tilde{\rho}$ such that $z_0\in\tilde{\Omega}$ and $u|_{\tilde{\Omega}}$ is maximal. Note that there is $\tilde{\varepsilon}>0$ such that a function $\rho-\tilde{\varepsilon}\tilde{\rho}$ is plurisubharmonic in some neighborhood of ${\rm cl}(\tilde{\Omega})$. Hence a function $\tilde{v}=\max\{v+\varepsilon\rho,v+\varepsilon(\rho-\tilde{\rho}\tilde{\varepsilon})\}$ is also plurisubharmonic and $\tilde{v}-u$ attains maximum only in some compact subset of $\tilde{\Omega}$ which is impossible because $u|_{\tilde{\Omega}}$ is maximal.
As $\varepsilon$ and $v$ were arbitrary we can conclude that the function $u$ is maximal. $\;\Box$

In \cite{h-l} the authors consider problem (\ref{MP}) for $\mathcal{F(J)}$-harmonic functions which they define in a different way than we define maximal functions but we will see that these concepts agree.

Let $\Omega\subset M$ and $u\in\mathcal{C\cap PSH}(\Omega)$. We say that $u$ is $\mathcal{F(J)}$-harmonic if for every $U\subset\Omega$ and for every smooth strictly plurisubharmonic function $\phi\leq u$ on $U$ we have $\phi<u$ on $U$. One can show (using the comparison principle) that $\mathcal{C}^2$ $\mathcal{F(J)}$-harmonic functions are exactly $\mathcal{C}^2$ solutions of (\ref{ma}) with $f=0$. 

\begin{prp} Let $\Omega$ and $u$ be as above.
Then\\
 i) If $u$ is maximal then $u$ is $\mathcal{F(J)}$-harmonic;
\\ ii) If $\Omega$ is such that there is a bounded strictly plurisubharmonic function $\rho\in\mathcal{C}^2(\Omega)$ and $u$ is $\mathcal{F(J)}$-harmonic then $u$ is maximal. 
\end{prp} 

{\it Proof:} The first assertion follows from definitions. To proof ii) we can assume, by Proposition \ref{lok}, that $\Omega$ is a smooth strictly pseudoconvex domain with defining function $\rho$  such that $u\in\mathcal{C}(\bar{\Omega})$. Let $\varepsilon>0$. By Theorem \ref{maksymalne} there is a continous maximal plurisubharmonic function $u_0$ equal on $\po$ to $u$. By Theorem \ref{mine theorem} there is a smooth strictly plurisubharmonic function $u_1$ such that $u-\varepsilon<u_1<u$ on a boundary and $$(i\partial\bar{\partial}u_1)^n=\frac{1}{2}\varepsilon^n(i\partial\bar{\partial}\rho)^n.$$ Then using the comparison principle (Proposition \ref{cp'}) we obtain $$u_0+\varepsilon\rho-\varepsilon\leq u_1\leq u\leq u_0$$ and thus we get $u=u_0$. $\;\Box$

$\newline$\textbf{Acknowledgments.} The author would like to express his
gratitude to Z. Błocki for helpful discussions and advice during
the work on this paper.


\begin{thebibliography}{9999999}
\bibitem[B]{b} F. Bertrand,  \textit{Local analysis on almost complex manifolds},
 PhD thesis 2007,\\ ( available on http://homepage.univie.ac.at/florian.bertrand/these.pdf ),
\bibitem[B-G]{b-g} F. Bertrand, H. Gaussier \textit{On the Gromov hyperbolicity of strongly pseudoconvex domains in almost complex manifolds}, arXiv:1202.4985,


\bibitem[C-K-N-S]{c-k-n-s}
   L. Caffarelli, J. J. Kohn, L. Nirenberg, J. Spruck
   \textit{The Dirichlet problem for non-linear second
order elliptic equations II: Complex Monge-Amp\`{e}re, and
uniformly elliptic equations}, Comm. Pure Appl. Math. 38 (1985),
209-252,
\bibitem[H-L]{h-l}
R. Harvey, B. Lawson 
\textit{Potential Theory on Almost Complex Manifolds}, ArXiv:1107.2584,
\bibitem[I-R]{i-r}
   S. Ivashkovich,  J.-P. Rosay,
   \textit{Schwarz-type lemmas for solutions of $\overline\partial$-inequalities and complete hyperbolicity of almost complex manifolds},
Ann. Inst. Fourier (Grenoble) 54 (2004), no. 7, 2387-2435,
\bibitem[M]{m} D. McDuff, \textit{Symplectic manifolds with contact type boundaries}, Invent. Math. 103 (1991), 651-671,

\bibitem[P]{p} N. Pali,  \textit{Fonctions plurisousharmoniques et courants positifs de type $(1, 1)$ sur une vari\'{e}t\'{e} presque
complexe},
Manuscripta Math. 118 (2005), no. 3, 311-337,
\bibitem[R]{r}
J.-P. Rosay
\textit{Notes on the Diederich-Sukhov-Tumanov normalization for almost complex structures}
Collect. Math. 60 (2009), no. 1, 43--62
   

 \end{thebibliography}
\end{document}